\documentclass{llncs}
\usepackage{amssymb, amsmath}
\usepackage{graphicx,color}
\newcommand{\boxi}{\ensuremath{box}}
\begin{document}
\date{}
\pagestyle{plain}
\title{Boxicity of Line Graphs}
\author{L. Sunil Chandran\inst{1} \and Rogers Mathew\inst{1} \and Naveen Sivadasan\inst{2}}
\institute{Department of Computer Science and Automation, \\ Indian Institute of
Science, \\
Bangalore - 560012, India. \\ \email{\{sunil,rogers\}@csa.iisc.ernet.in} \and 
Department of Computer Science and Engineering, \\
Indian Institute of Technology, \\
Hyderabad - 502205, India. \\ \email{nsivadasan@iith.ac.in}}
\maketitle
\bibliographystyle{plain}
\begin{abstract}
Boxicity of a graph $H$, denoted by $\boxi(H)$, is the minimum integer $k$ such that 
$H$ is an intersection graph of axis-parallel $k$-dimensional 
boxes in $\mathbb{R}^k$. 
In this paper, we show that for a line graph $G$ of a multigraph, 
$\boxi(G) \leq 2\Delta(\lceil \log_2 \log_2 \Delta \rceil + 3) + 1$, 
where $\Delta$ denotes the maximum degree of $G$.  
Since $\Delta \leq 2(\chi - 1)$, for any line graph $G$ with chromatic 
number $\chi$, $\boxi(G) = O(\chi \log_2 \log_2 (\chi) )$. 
For the $d$-dimensional hypercube $H_d$, we prove that 
$ \boxi(H_d) \geq \frac{1}{2}\left(\lceil \log_2 \log_2 d \rceil + 1\right)$. 
The question of finding a non-trivial lower bound for $\boxi(H_d)$ was left open 
by Chandran and Sivadasan in [L. Sunil Chandran and Naveen Sivadasan. 
The cubicity of Hypercube Graphs. 
\emph{Discrete Mathematics}, 308(23):5795-5800, 2008]. 

The above results are consequences of bounds 
that we obtain for the boxicity of fully subdivided graphs
(a graph which can be obtained by subdividing every edge of a graph exactly once). 

\end{abstract}
\section{Introduction}
Given a family $\mathcal{F}$ of sets, a graph $G=(V,E)$ is called an 
\emph{intersection graph} of sets from $\mathcal{F}$, if there exists a 
map $f~:~V(G) \rightarrow \mathcal{F}$ such that 
$(u,v) \in E(G) \Leftrightarrow f(u) \cap f(v) \neq \emptyset$. If the sets 
in $\mathcal{F}$ are intervals on a real line, then we call $G$ an \emph{interval 
graph}. In other words, interval graphs are intersection graphs of 
intervals on the real line. In $\mathbb{R}^k$, an axis parallel $k$-dimensional 
box or a \emph{$k$-box} is a cartesian product $R_1 \times R_2 \times \cdots \times R_k$, where 
each $R_i$ is a closed interval $[a_i,b_i]$ on the real line. 
A graph $G$ is said to have a $k$-box representation if there exists a mapping from the  
vertices of $G$ to $k$-boxes in the $k$-dimensional eucledian space 
such that two vertices in $G$ are adjacent if and only if their corresponding $k$-boxes  
have a non-empty intersection. 
\emph{Boxicity} of $G$, denoted 
by $\boxi(G)$, is 
the minimum positive integer $k$ such that $G$ has a 
$k$-box representation. As each interval can also be viewed as an 
axis parallel $1$-dimensional box, interval graphs are precisely the 
class of graphs with boxicity 1. 
We take the boxicity of a complete graph to be 1. 
\subsection{Background}
The concept of boxicity was introduced by F.S. Roberts in 1969 \cite{Roberts}.
Cozzens \cite{Coz} showed that computing the boxicity of a graph is 
NP-hard. Yannakakis in \cite{Yan1} improved this result. Finally, 
Kratochvil \cite{Kratochvil} showed that deciding whether the boxicity 
of a graph is at most 2 itself is NP-complete.

Box representation of graphs finds application in niche overlap (competition) in
ecology and to problems of fleet maintenance in operations research
(see \cite{CozRob}). Given a low dimensional box representation,  
some well known NP-hard problems become polynomial time solvable. 
For instance, the max-clique problem is polynomial time solvable for graphs with 
boxicity $k$ because the number of maximal cliques in such graphs is only 
$O((2n)^k)$.  
 
Roberts proved that for every graph $G$ 
on $n$ vertices, $\boxi(G) \leq \lfloor \frac{n}{2} \rfloor$. He gave a tight example to this  
by showing that a complete $\frac{n}{2}$-partite graph with 2 vertices in each part 
has its boxicity equal to $\frac{n}{2}$. In \cite{chintan}, it was shown that if   
$t$ denotes the size of a 
minimum vertex cover of $G$, then $\boxi(G) \leq \lfloor \frac{t}{2}\rfloor + 1$. 
Chandran, Francis and Sivadasan showed in \cite{tech-rep} that, 
for any graph $G$ on $n$ vertices having maximum degree $\Delta$,
$\boxi(G) \leq (\Delta + 2)\ln n $. An upper bound solely in terms of the maximum
degree $\Delta$, which says $\boxi(G) \leq 2\Delta^2$, is proved in \cite{CFNMaxdeg}.
Esperet \cite{Esperet} improved this bound to $\Delta^2 + 2$. Recently 
Adiga, Bhowmick and Chandran \cite{DiptAdiga} showed that $\boxi(G)=O(\Delta \log^2 \Delta)$. 
Chandran and Sivadasan in \cite{CN05} found a relation between treewidth and boxicity 
which says $\boxi(G) \leq \mathrm{tw}(G)+ 2$, 
where $\mathrm{tw}(G)$ denotes the treewidth of graph $G$.

Attempts on finding better bounds for boxictiy of special graph 
classes can also be seen in the literature. Scheinerman \cite{Scheiner} 
showed that outerplanar graphs have boxicity at most 2. Thomassen \cite{Thoma1} 
proved that the boxicity of planar graphs is not greater than 3. Cozzens 
and Roberts \cite{CozRob} have done a study on the boxicity of split graphs. 
Results on the boxicity of Chordal graphs, AT-free graphs, permutation graphs 
etc. can be seen in \cite{CN05}. Better bounds for   
the boxicity of Circular Arc graphs and AT-free graphs can be seen in \cite{Dipt,Bhowmick20101536}. 
In \cite{SunMatRog} it was shown that, there exist chordal bipartite graphs with  
arbitrarily high boxicity.  
\subsection{An Equivalent Definition for Boxicity}
Let $G,G_1,G_2,\ldots,G_b$ be a collection of graphs with $V(G) = V(G_i)$,  
for any $i \leq b$. We say $G= \bigcap_{i=1}^{b}G_i$ 
when $E(G) = \bigcap_{i=1}^{b}E(G_i)$. 
The following lemma gives the relationship between interval graphs 
and intersection graphs of $k$-boxes.
\begin{lemma}[Roberts\cite{Roberts}]
\label{Robertslemma}
For any graph $G$, $\boxi(G) \leq k$ if and only if
there exist $k$ interval graphs $I_1, I_2, \ldots, I_k$
 such that $G = \bigcap_{i=1}^{k}I_i$.
\end{lemma}
From the above lemma, we can say that boxicity of a graph $G$ is the minimum positive integer $k$ for which 
there exist $k$ interval graphs $I_1,I_2 \ldots, I_k$ such that $G= \bigcap_{i=1}^kI_i$.

We have seen that intervals graphs are intersection graphs of 
intervals on the real line. Hence for any interval graph 
$I$, there exists a map $f~:~V(I)\rightarrow \{X \subseteq \mathbb{R}~|~
X\mbox{ is a closed interval}\}$ such that, for any $u,v \in V(I)$, 
$(u,v) \in E(I)$ if and only if $f(u) \cap f(v) \neq \emptyset$. 
Such a map $f$ is called an \emph{interval representation} of $I$. 
An interval graph can have more than one interval representation. 
It is known that given an interval graph $I$, we can find an interval 
representation for $I$ in which no two intervals share any endpoints.  
\subsection{Preliminaries}
Except in Theorem \ref{linegraphtheorem}, Section \ref{linegraphsection}, we consider only finite, undirected,  
and simple graphs. In Theorem \ref{linegraphtheorem}, we consider finite, undirected 
multigraphs. For any finite positive integer $n$, let $[n]$ denote the 
set $\{1,2,\ldots n\}$. 
For a graph $G$, we use $V(G)$ and $E(G)$ to denote the set of 
its vertices and edges respectively. For any $v \in V(G)$, $N_G(v) := \{u~|~(v,u) \in E(G)\}$
and $d_G(v) := |N_G(v)|$. The maximum degree of $G$ is denoted by $\Delta(G)$. 
$\chi(G)$ represents the chromatic number of $G$. 
We say that an edge $e_i$ is a neighbour of 
another edge $e_j$ in $G$, if they share an endpoint. 
Given two graphs $G$ and $H$, we say $G=H$ when $G$ is isomorphic to $H$.  

We say that a graph $G$ is obtained by \emph{fully subdividing} $H$, if $G$ is  
obtained as a result of subdividing every edge of $H$ exactly once.
Given a multigraph $H$, we define a graph $L(H)$ in the 
following way: $V(L(H)) = E(H)$ and $E(L(H)) = 
\{(e_1,e_2)~|~e_1,e_2 \in E(H),~e_1 \mbox{ and } e_2 \mbox{ share an  
endpoint in }$ $ H \}$. A graph $G$ is a \emph{line graph} if and 
only if there exists a multigraph $H$ such that $G$ is isomorphic to $L(H)$.
Let $I$ be an interval graph and $f$ an interval representation of $I$. 
Then, $\forall x\in V(I)$,  we use $l(f(x))$ and $r(f(x))$ to denote 
the left and right endpoint respectively of the interval $f(x)$. 
\subsection{Our Results}
In this paper, we show that for a line graph $G$ with maximum degree 
$\Delta$,
$$\boxi(G) \leq 2\Delta(\lceil \log_2 \log_2 \Delta \rceil + 3) + 1.$$
From the above result, we also infer that if 
chromatic number of $G$ is $\chi$, then $\boxi(G) = O(\chi \log_2 \log_2 (\chi))$. 
Recall that, in \cite{DiptAdiga} it was shown that for any graph $G$, 
$\boxi(G)\leq c\cdot \Delta \log^2 \Delta$, where $c$ is a large constant. 
Hence, for the class of line graphs, 
our result is an improvement over the best bound known for general graphs. 
Moreover, in contrast with the result in \cite{DiptAdiga}, the proof here 
is constructive and easily gives an efficient algorithm to get a  
box representation for the given line graph. We leave the tightness of our result open.

The main 
supporting result that we have used to prove the above result is the 
following (this itself may be independently interesting): For a graph $G$ 
obtained by fully subdividing another graph $H$, 
$\boxi(G) \leq \lceil \log_2 \log_2 (\Delta) \rceil + 3$, where 
$\Delta$ is the maximum degree of $G$. At the end of the paper, we point 
out another consequence of this supporting result. 
For the $d$-dimensional hypercube $H_d$, 
$$ \boxi(H_d) \geq \frac{\lceil \log_2 \log_2 d \rceil + 1}{2}.$$
It was shown by Chandran and Sivadasan in \cite{CN99} that 
$\boxi(H_d) \leq \frac{cd}{\log d}$, where $c$ is a constant. They had raised 
the question of finding a non-trivial lower bound for $\boxi(H_d)$.
\section{Boxicity of a Fully Subdivided Complete Graph}
Let $S=\{\sigma_1, \sigma_2, \ldots , \sigma_p\}$ be a set of permutations of 
$[n]$, where $n$ is any finite positive integer. 
$S$ is called \emph{$k$-suitable} for $[n]$ if for any $k$-element subset 
$X \subseteq [n]$ and for any $x \in X$, there exists a permutation $\sigma \in S$ with the 
following property: 
\[\sigma^{-1}(x) \geq \sigma^{-1}(y), \forall y \in X.\] 
The minimum cardinality of a $k$-suitable set for $[n]$ is 
denoted by $N'(n,k)$. Spencer \cite{scramble} proved that 
\[N'(n,3) < \log_2 \log_2 n + \frac{1}{2} \log_2 \log_2 \log_2 n + \log_2 (\sqrt{2} \pi) + o(1).\]
In this paper, we are interested in a slightly relaxed version of the notion 
of $3$-suitability. Given a permutation $\sigma$ of $[n]$ and 
$s,t \in [n]$, let 
\begin{eqnarray}
\label{betaequation}
\beta(s,t,\sigma) & = & \{ x~|~\sigma^{-1}(s) < \sigma^{-1}(x) < \sigma^{-1}(t) \nonumber \\
 &  & \mbox{ or } \sigma^{-1}(t) < \sigma^{-1}(x) < \sigma^{-1}(s) \}. 
\end{eqnarray}
A set $S=\{\sigma_1, \sigma_2, \ldots , \sigma_p\}$ is called \emph{simply $3$-suitable}
for $[n]$, if for each pair $s,t \in [n]$, $\bigcap_{i=1}^p \beta(s,t,\sigma_i) = \emptyset$. 
In other words, for every triple $x,s,t \in [n]$ there exists a permutation $\sigma \in S$ 
such that either $\sigma^{-1}(x) < \min \left(\sigma^{-1}(s), \sigma^{-1}(t)\right)$ or 
$\sigma^{-1}(x) > \max \left(\sigma^{-1}(s), \sigma^{-1}(t)\right)$. It 
is easy to see that any $3$-suitable set is also a simply $3$-suitable set while the 
converse is clearly not true. Let $N(n)$ be the minimum possible cardinality of a 
simply $3$-suitable set for $[n]$. From Spencer's bound on 
$N'(n,3)$, we have $N(n) \leq N'(n,3) < 
\log_2 \log_2 n + \frac{1}{2} \log_2 \log_2 \log_2 n + \log_2 (\sqrt{2} \pi) + o(1)$. 
But since simply $3$-suitability is a more relaxed notion than $3$-suitability, 
we can get the following exact formula for $N(n)$:
\begin{lemma}
\label{scramblelemma}
$N(n) = \lceil \log_2 \log_2 n \rceil + 1$. 
\end{lemma}
\begin{proof}
Erd\H{o}s and Szekeres \cite{ErdosSzekeres} proved that if $\sigma_1$ and 
$\sigma_2$ are two permutations of $[n^2 + 1]$, then there exists 
some $X \subset [n^2 + 1]$ with $|X|=n+1$ such that the permutation 
of $X$ obtained by restricting $\sigma_1$ to $X$ is the same as 
the permutation obtained by restricting $\sigma_2$ to $X$. By an 
easy inductive argument (as Spencer points out in \cite{scramble}) 
we can show that if $\sigma_1, \sigma_2, \ldots 
\sigma_{s+1}$ are permutations of $[2^{2^s} + 1]$, then there exists some 
triple $\{x,y,z\}$ such that the order of these 3 elements with respect 
to each permutation $\sigma_1, \sigma_2, \ldots \sigma_{s+1}$ is the same.  
This implies that $N(n) \geq\lceil \log_2 \log_2 n \rceil + 1$. 

We need to show that when $n\leq 2^{2^i}$, $N(n) \leq i+1$.
Note that when the permutations in a simply 3-suitable set $S$ for $[n]$ are restricted to $[n_1]$ 
(where $n_1 < n$), $S$ becomes a simply 3-suitable set for $[n_1]$. Hence it is 
enough to prove that, when $n=2^{2^i}$, $N(n) \leq i+1$.
We prove this by induction on $i$. The base case, when $i=0$ and $n=2$, is
 trivially true. For any $i<i_1$, assume $N(n) \leq i+1$. 
Let $i=i_1$, $n=2^{2^{i_1}}$ and $n_1 = 2^{2^{i_1-1}}$. Then 
$n=n_1 \cdot n_1$. So set $[n]$ can be partitioned into $n_1$ sets 
$A_1, A_2, \ldots A_{n_1}$, where for any $p \in [n_1]$, 
$A_p = \{(p-1)n_1 + 1, (p-1)n_1 + 2, \ldots , (p-1)n_1 + n_1\}$. 
Clearly for any $a \in [n]$, there exist $k,p \in [n_1]$  
such that $a=(p-1)n_1 + k$.
By induction hypothesis, there exists a simply 3-suitable set 
$S'= \{\eta_1, \eta_2, \ldots \eta_{i_1}\}$  of $[n_1]$. 
Then we define $i_1 +1$ permutations 
$S=\{\sigma_1, \ldots, \sigma_{i_1+1}\}$ for $[n]$ as follows:
\begin{eqnarray*}
\label{generatescramble1}
\sigma_j^{-1}(a) & = & (\eta_j^{-1}(p) - 1)n_1 + \eta_j^{-1}(k) \mbox{, where } 1 \leq j \leq i_1. \\
\label{generatescramble2}
\sigma_{i_1+1}^{-1}(a) & = & (n_1 - \eta_{i_1}^{-1}(p))n_1 + \eta_{i_1}^{-1}(k).
\end{eqnarray*}

We claim that $S$ is a simply 3-suitable set for $[n]$ i.e., for any
$s,t \in [n]$, $\bigcap_{i=1}^{i_1+1}\beta(s,t,\sigma_i)=\emptyset$. 
Let $s \in A_p$ and $t \in A_q$. Consider the 2 cases below: \\
\textbf{case 1}: If $p=q$, then there exist $k_1,k_2 \in [n_1]$  
with $k_1 \neq k_2$ such that, $s=(p-1)n_1+ k_1$ and $t=(p-1)n_1 + k_2$. 
Consider a permutation $\sigma_j$, where $j \in [i_1]$. 
\begin{eqnarray*}
\beta(s,t,\sigma_j) & = & \{ x~|~\sigma_j^{-1}(s) < \sigma_j^{-1}(x) < \sigma_j^{-1}(t) \\
 &  & \mbox{ or } \sigma_j^{-1}(t) < \sigma_j^{-1}(x) < \sigma_j^{-1}(s) \} \\
 & = & \{ x~|~(\eta_j^{-1}(p)-1)n_1 + \eta_j^{-1}(k_1) < \sigma_j^{-1}(x) < 
(\eta_j^{-1}(p)-1)n_1 + \eta_j^{-1}(k_2) \\
 &  & \mbox{ or } (\eta_j^{-1}(p)-1)n_1 + \eta_j^{-1}(k_2) < \sigma_j^{-1}(x) < 
(\eta_j^{-1}(p)-1)n_1 + \eta_j^{-1}(k_1) \}. 
\end{eqnarray*}
If $\beta(s,t,\sigma_j) \neq \emptyset$, then consider any $x \in \beta(s,t,\sigma_j)$. 
Clearly $x \in A_p$. Let $x = (p-1)n_1 + k_3$. 
From the above, it is clear that either $\eta_j^{-1}(k_1) < \eta_j^{-1}(k_3) < \eta_j^{-1}(k_2)$ 
or $\eta_j^{-1}(k_2) < \eta_j^{-1}(k_3) < \eta_j^{-1}(k_1)$. This means that 
$x \in \beta(s,t,\sigma_j) \implies k_3 \in \beta(k_1,k_2,\eta_j)$. Therefore, 
$\bigcap_{j=1}^{i_1}\beta(s,t,\sigma_j) \neq \emptyset \implies 
\bigcap_{j=1}^{i_1}\beta(k_1,k_2,\eta_j) \neq \emptyset$. 
By induction hypothesis, we know that $\bigcap_{j=1}^{i_1}\beta(k_1,k_2,\eta_j) = 
\emptyset$. Hence $\bigcap_{j=1}^{i_1}\beta(s,t,\sigma_j) = \emptyset$.    \\
\textbf{case 2}: If $p \neq q$, then $\exists k_1,k_2 \in [n_1]$ such that 
$s=(p-1)n_1 + k_1$ and $t=(q-1)n_1 + k_2$. Let $x = (r-1)n_1 + k_3$. 
Now $x \in \bigcap_{j=1}^{i_1}\beta(s,t,\sigma_j)$ implies, for any 
$j \in [n_1]$, 
$
 (\eta_j^{-1}(p)-1)n_1 + \eta_j^{-1}(k_1)  <  (\eta_j^{-1}(r)-1)n_1 + \eta_j^{-1}(k_3)  <   
(\eta_j^{-1}(q)-1)n_1 + \eta_j^{-1}(k_2) 
 \mbox{ or } 
(\eta_j^{-1}(q)-1)n_1 + \eta_j^{-1}(k_2)  <  (\eta_j^{-1}(r)-1)n_1 + \eta_j^{-1}(k_3)  <   
(\eta_j^{-1}(p)-1)n_1 + \eta_j^{-1}(k_1)
$. It follows that $\eta_j^{-1}(p) \leq \eta_j^{-1}(r) \leq \eta_j^{-1}(q)$ or 
$\eta_j^{-1}(q) \leq \eta_j^{-1}(r) \leq \eta_j^{-1}(p)$. If $r \notin \{p,q\}$,  
then $\eta_j^{-1}(p) < \eta_j^{-1}(r) < \eta_j^{-1}(q)$ or 
$\eta_j^{-1}(q) < \eta_j^{-1}(r) < \eta_j^{-1}(p)$ i.e., 
$ r \in \bigcap_{j=1}^{i_1}\beta(p,q,\eta_j)$ which contradicts the induction 
hypothesis that $\bigcap_{j=1}^{i_1}\beta(p,q,\eta_j)= \emptyset$. 

Therefore we infer that $r=p$ or $r=q$. Let 
$r=p$ (proof is similar when $r=q$). If $x \in \bigcap_{j=1}^{i_1+1}\beta(s,t,\sigma_j)$ 
then we have $x \in \beta(s,t,\sigma_{i_1})$ and therefore  
$\sigma_{i_1}^{-1}(s) < \sigma_{i_1}^{-1}(x) < \sigma_{i_1}^{-1}(t)$ or 
$\sigma_{i_1}^{-1}(t) < \sigma_{i_1}^{-1}(x) < \sigma_{i_1}^{-1}(s)$. 
Without loss of generality, let 
$\sigma_{i_1}^{-1}(s) < \sigma_{i_1}^{-1}(x) < \sigma_{i_1}^{-1}(t)$. Then 
$(\eta_{i_1}^{-1}(p)-1)n_1 + \eta_{i_1}^{-1}(k_1)  <  (\eta_{i_1}^{-1}(r)-1)n_1 + \eta_{i_1}^{-1}(k_3)  <   
(\eta_{i_1}^{-1}(q)-1)n_1 + \eta_{i_1}^{-1}(k_2)$. Since $p=r$, we have 
$\eta_{i_1}^{-1}(p) = \eta_{i_1}^{-1}(r)$ and therefore $\eta_{i_1}^{-1}(k_1) < 
\eta_{i_1}^{-1}(k_3)$. This also allows us to infer that 
$(n_1 - \eta_{i_1}^{-1}(p))n_1 + \eta_{i_1}^{-1}(k_1) < 
(n_1 - \eta_{i_1}^{-1}(r))n_1 + \eta_{i_1}^{-1}(k_3)$. That is 
$\sigma_{i_1+1}^{-1}(s) <  \sigma_{i_1+1}^{-1}(x)$. On the other hand, 
$(n_1 - \eta_{i_1}^{-1}(q))n_1 + \eta_{i_1}^{-1}(k_2) < 
(n_1 - \eta_{i_1}^{-1}(p))n_1 + \eta_{i_1}^{-1}(k_1)$ (since 
$\eta_{i_1}^{-1}(p) < \eta_{i_1}^{-1}(q)$). Therefore, 
$\sigma_{i_1+1}^{-1}(t) <  \sigma_{i_1+1}^{-1}(s)$. So we have, 
$\sigma_{i_1+1}^{-1}(t) < \sigma_{i_1+1}^{-1}(s) <  \sigma_{i_1+1}^{-1}(x)$. Hence 
$x \notin \beta(s,t,\sigma_{i_1+1})$ contradicting our assumption that 
$x \in \bigcap_{j=1}^{i_1+1}\beta(s,t,\sigma_j)$.
\qed
\end{proof} 
\begin{theorem}
\label{completegraphtheorem}
Let $G$ be the graph obtained by fully subdividing the complete graph $K_n$. Then 
$\frac{\lceil \log_2 \log_2 n \rceil + 1}{2} \leq \boxi(G) \leq \lceil \log_2 \log_2 n \rceil + 2$. 
\end{theorem}
\begin{proof}
Let $v_1, v_2, \ldots v_n$ be the vertices of $K_n$ and $
e_1, e_2, \ldots e_m$ its edges, where $m =$  $n\choose2$. Let $u_{p\cdot q}$ 
denote the vertex introduced when subdividing the edge $(v_p,v_q) \in E(K_n)$, where 
$p<q$. Thus the graph $G$ obtained by fully subdividing $K_n$ has the vertex set 
$V(G) = \{v_1,v_2,\ldots v_n\} \cup \{u_{p\cdot q}~|~1 \leq p < q \leq n\}$ 
and $E(G) = \{(v_p,u_{p\cdot q})~|~1 \leq p < q \leq n\} \cup \{(v_q,u_{p\cdot q}~|~ 
1 \leq p<q\leq n)\}$. 

We first show that $\boxi(G) \leq \lceil \log_2 \log_2 n \rceil + 2$. 
Let $k=\lceil \log_2 \log_2 n \rceil + 1$. By Lemma \ref{scramblelemma},  
there exists a simply 3-suitable set $S=\{\sigma_1, \ldots, \sigma_k\}$ for 
$[n]$. Using $S$, we construct a $(k+1)$-dimensional 
box representation for $G$. Corresponding to each permutation $\sigma_i$ of 
$[n]$ in $S$, we construct an interval graph $I_i$ as follows. Let $f_i$ denote the interval 
representation of $I_i$. 
\begin{eqnarray*}
\mbox{for every }v_p \in V(G), & f_i(v_p) & =  [\sigma_i^{-1}(p),\sigma_i^{-1}(p)]. \\
\mbox{for every }u_{p\cdot q} \in V(G), & f_i(u_{p\cdot q}) & =  
[\sigma_i^{-1}(p),\sigma_i^{-1}(q)] \mbox{, if } \sigma_i^{-1}(p)] < \sigma_i^{-1}(q). \\
\mbox{for every }u_{p\cdot q} \in V(G), & f_i(u_{p\cdot q}) & = 
[\sigma_i^{-1}(q),\sigma_i^{-1}(p)] \mbox{, if } \sigma_i^{-1}(q) < \sigma_i^{-1}(p)]. 
\end{eqnarray*}
The interval representation $f_{k+1}$ of the $(k+1)$th interval graph $I_{k+1}$ is as follows: 
\begin{eqnarray*}
\mbox{for every }v_p \in V(G), & f_{k+1}(v_p) & = [1,m]. \\
\mbox{for every }u_{p\cdot q} \in V(G), & f_{k+1}(u_{p\cdot q}) & = [j,j] \mbox{, where } u_{p\cdot q} \mbox { 
was obtained by } \\
& & \mbox{subdividing edge } e_j = (v_p,v_q) \mbox{ of } K_n.  
\end{eqnarray*}
By Lemma \ref{Robertslemma}, in order to prove that $\boxi(G) \leq k+1$ it is sufficient to 
show that $\bigcap_{i=1}^{k+1}I_i = G$, i.e.,\\
(i) each $I_j$ is a supergraph of $G$. \\
(ii) for any $(x,y) \notin E(G)$, there exists some interval graph $I_i$ 
such that $(x,y) \notin E(I_i)$.

Recall that any edge of $G$ is of the form $(v_p,u_{pq})$ or 
$(v_q,u_{pq})$, where $v_p, v_q \in V(K_n)$. It is easy to verify 
that, for any $i \in [k+1]$, $f_i(u_{pq}) \cap f_i(v_p) \neq \emptyset$ 
and $f_i(u_{pq}) \cap f_i(v_q) \neq \emptyset$. Therefore (i) is true. 

Let $(x,y) \notin E(G)$. In order to prove (ii), we consider the following cases: 
%
\textbf{case 1}: $x=v_p,y= v_q$, for some $1\leq p< q \leq n$. \\
It is easy to see that $f_1(v_p) \cap f_1(v_q) = \emptyset$ and therefore 
$(v_p,v_q) \notin E(I_1)$. \\
\textbf{case 2}: $x=u_{p\cdot q}, y = u_{r\cdot s}$  and $u_{p\cdot q} \neq u_{r\cdot s}$. \\
Clearly, $f_{k+1}(u_{p\cdot q}) \cap f_{k+1}(u_{r\cdot s}) = \emptyset$ and therefore 
$(u_{p\cdot q}, u_{r\cdot s}) \notin E(I_{k+1})$. \\
\textbf{case 3}: $x=v_p, y=u_{r\cdot s}$, for any $p,r,s \in [n]$, 
$p \notin \{r,s\}$ and $r < s$. \\
Since $S$ is a simply 3-suitable set for $[n]$ there exists 
a permutation $\sigma_j$ such that $p \notin \beta(r,s,\sigma_j)$ i.e., 
either $\sigma_j^{-1}(p) < \min(\sigma_j^{-1}(r),\sigma_j^{-1}(s))$ or 
$\sigma_j^{-1}(p) > \max(\sigma_j^{-1}(r),\sigma_j^{-1}(s))$. 
Now it is easy to see that, $f_j(v_p) \cap f_j(u_{r\cdot s}) = 
\emptyset$ and therefore $(v_p, u_{r\cdot s}) \notin E(I_j)$.  
We thus prove (ii) and thereby 
prove that $\boxi(G) \leq \lceil \log_2\log_2n \rceil + 2$.


We now show that $\boxi(G) \geq \frac{\lceil \log_2 \log_2 n \rceil + 1}{2}$. Let 
$\boxi(G) = b$. By Lemma \ref{Robertslemma} there  
exist $b$ interval graphs, say $I_1, I_2, \ldots , I_b$, such that 
$G = \bigcap_{i=1}^{b}I_i$. For any $i \in [b]$, let $f_i$ be an 
interval representation of $I_i$ such that no two intervals share any endpoints. 
From each $f_i$, generate two permutations $L_i$ and $R_i$ of $[n]$ in the following way. 
For $p,q \in [n]$, $p\neq q$, $L_i^{-1}(p) < L_i^{-1}(q) \Leftrightarrow l(f_i(v_p)) < 
l(f_i(v_q))$. Similarly, $R_i^{-1}(p) < R_i^{-1}(q) \Leftrightarrow r(f_i(v_p)) < 
r(f_i(v_q))$

Consider the set $S=\{L_1,R_1,L_2,R_2,\ldots L_b,R_b\}$ of permutations 
of $[n]$. We claim that $S$ is a simply 3-suitable set for $[n]$.   
Let $s,t \in [n]$. Then for any $i \in [b]$, 
%
\begin{eqnarray}
\label{implications1}
x \in \beta(s,t,L_i) & \implies & 
\left(L_i^{-1}(s) < L_i^{-1}(x) < L_i^{-1}(t)\right) \mbox{ or} \\
&  & \left(L_i^{-1}(t) < L_i^{-1}(x) < L_i^{-1}(s)\right) \nonumber \\ 
& \implies & \left(l(f_i(v_s)) < l(f_i(v_x)) < l(f_i(v_t))\right) \mbox{ or } \nonumber \\
&  & \left(l(f_i(v_t)) < l(f_i(v_x)) < l(f_i(v_s))\right). \nonumber \\
\label{implications2}
x \in \beta(s,t,R_{i}) & \implies & 
\left(R_{i}^{-1}(s) < R_{i}^{-1}(x) < R_{i}^{-1}(t)\right) \mbox{ or} \\
&  & \left(R_{i}^{-1}(t) < R_{i}^{-1}(x) < R_{i}^{-1}(s)\right) \nonumber \\ 
& \implies & \left(r(f_i(v_s)) < r(f_i(v_x)) < r(f_i(v_t))\right) \mbox{ or } \nonumber \\
&  & \left(r(f_i(v_t)) < r(f_i(v_x)) < r(f_i(v_s))\right). \nonumber 
\end{eqnarray}
Suppose, for contradiction, $x \in \bigcap_{j=1}^{b}\left(\beta(s,t,L_j)
\cap \beta(s,t,R_j)\right)$. Consider any $i \in [b]$. 
Let $y=\max(l(f_i(v_s)),l(f_i(v_t)))$ and 
$z=\min(r(f_i(v_s)),r(f_i(v_t)))$. 
Consider the two cases below: \\
\textbf{case 1}: $y<z$. 
Then by implications (\ref{implications1}) and (\ref{implications2}) it is clear that $l(f_i(v_x)) < y = 
\max(l(f_i(v_s)),l(f_i(v_t)))$ and $r(f_i(v_x)) > 
z=\min(r(f_i(v_s)),r(f_i(v_t)))$. Therefore, $[y,z] \subseteq f_i(v_x)$. 
Now we will show that $f_i(u_{s\cdot t}) \cap [y,z] \neq \emptyset$ which will 
immediately imply that $f_i(u_{s\cdot t}) \cap f_i(v_x) \neq \emptyset$. 
If $f_i(u_{s\cdot t}) \cap [y,z] = \emptyset$, then either 
$r(f_i(u_{s\cdot t})) < y$ or $l(f_i(u_{s\cdot t})) > z$. In both these cases, 
it is easy to see that either $(u_{s\cdot t}, v_s) \notin E(I_i)$ or 
$(u_{s\cdot t}, v_t) \notin E(I_i)$. This contradicts the fact that $I_i$ is a 
supergraph of $G$. Hence $f_i(u_{s\cdot t}) \cap [y,z] \neq \emptyset$ and therefore
$(u_{s\cdot t}, v_x) \in E(I_i)$.\\ 
\textbf{case 2}: $y>z$. Since $(u_{s\cdot t}, v_s) \in E(I_i)$ and 
$(u_{s\cdot t}, v_t) \in E(I_i)$, we have $r(f_i(u_{s\cdot t})) > y$ and 
$l(f_i(u_{s\cdot t})) < z$. Therefore, $[z,y] \subseteq f_i(u_{s\cdot t})$. 
Now we will show that $f_i(v_x) \cap [z,y] \neq \emptyset$ which will 
immediately imply that $f_i(u_{s\cdot t}) \cap f_i(v_x) \neq \emptyset$. 
If $f_i(v_x) \cap [z,y] = \emptyset$, then either 
$r(f_i(v_x)) < z$ or $l(f_i(v_x))) > y$. In both these cases, 
we contradict implications (\ref{implications1}) and (\ref{implications2}) which state that
 $r(f_i(v_x))$ is sandwiched between $r(f_i(v_s))$ and 
$r(f_i(v_t))$, and $l(f_i(v_x))$ is sandwiched between $l(f_i(v_s))$ and 
$l(f_i(v_t))$. Hence $f_i(v_x) \cap [z,y] \neq \emptyset$ and therefore 
$(u_{s\cdot t}, v_x) \in E(I_i)$.

Thus we conclude that if there exists an $x \notin \{s,t\}$ such that 
$x \in \bigcap_{j=1}^{2b}\beta(s,t,$ $\sigma_j)$, then 
$(u_{s\cdot t},v_x) \in E(\bigcap_{i=1}^b I_i)$  which implies that $(u_{s\cdot t},v_x) \in E(G)$. 
But this contradicts the fact that $(u_{s\cdot t},v_x) \notin E(G)$ and hence 
$\bigcap_{j=1}^{2b}\beta(s,t,\sigma_j) = \emptyset$ i.e., $S$ is a simply 3-suitable set. 
Then by Lemma \ref{scramblelemma}, $|S| = 2b \geq \lceil \log_2 \log_2 n \rceil + 1$ or 
 $\boxi(G) \geq \frac{\lceil \log_2 \log_2 n \rceil + 1}{2}$. 
\qed
\end{proof}
\remark{Louis Esperet informed us that he had independently observed Theorem 
\ref{completegraphtheorem}. But he has not published it. 
We thank him for personal communication. In \cite{Esperet}, he also conjectures that 
for any graph $G$, (i) $\boxi(G) \leq a(G) + \kappa$, (ii) $\boxi(G) \leq \lambda \cdot a(G)$, 
where $\kappa$, $\lambda$ are constants and $a(G)$ refers to the arboricity of $G$. As 
arboricity of any graph is upper bounded by its degeneracy and since fully subdivided 
complete graphs are 2-degenerate, Theorem \ref{completegraphtheorem} disproves 
Esperet's both conjectures.  
}
\section{Boxicity of a Fully Subdivided Graph of Chromatic Number $\chi$}
\begin{theorem}
\label{sparsegraphtheorem}
Let $H$ be a graph with chromatic number $\chi$ and  
let $G$ be the graph obtained by fully subdividing $H$. Then,
 $\boxi(G) \leq \lceil \log_2 \log_2 \chi \rceil + 3$. 
\end{theorem}
\begin{proof}
Given a colouring of $H$ using $\chi$ colours, let $C_1, C_2 \ldots C_\chi$ represent the 
$\chi$ colour classes. Let $|C_i| = c_i$ and $c_{max}=\max_i(c_i)$. 
Give an arbitrary order to the vertices in each colour class. Let $v_{ij}$ denote the 
$j$-th vertex in the $i$-th colour class, where $i \in [\chi]$ and $j \in [c_i]$. 
Let $E(H)=\{e_1, e_2, \ldots ,e_m\}$ be the edge set of $H$. Let 
$u_{pq\cdot rs}$ denote the vertex introduced while subdividing the edge $(v_{pq}, v_{rs})$,
 where $p < r$. 
Let $k = \lceil \log_2 \log_2 \chi \rceil + 1$. By Lemma \ref{scramblelemma}, there exists 
a simply 3-suitable set $S=\{\sigma_1, \ldots \sigma_k\}$ for $[\chi]$. 
We use $S$ to construct a $(k+2)$-dimensional box representation for $G$. 
Corresponding to each 
permutation $\sigma_i \in S$, we construct an interval graph $I_i$ as follows. 
Let $f_i$ denote the interval representation of $I_i$. \\ When $i \leq k$, 
\begin{eqnarray*}
\mbox{for every } v_{pq} \in E(G), & f_i(v_{pq}) & = [g_i(p,q),g_i(p,q)], \\
& & \mbox{where } g_i(p,q)=\sigma_i^{-1}(p) + \frac{q-1}{c_{max}}. \\
\mbox{for every } u_{pq\cdot rs} \in E(G), & f_i(u_{pq\cdot rs}) & = [g_i(p,q),g_i(r,s)] \mbox{, if } g_i(p,q) < g_i(r,s). \\
\mbox{for every } u_{pq\cdot rs} \in E(G), & f_i(u_{pq\cdot rs}) & = [g_i(r,s),g_i(p,q)] \mbox{, if } g_i(r,s) < g_i(p,q), \\ 
& & \mbox{where }g_i(p,q) = \sigma_i^{-1}(p) + \frac{q-1}{c_{max}} \\
& & \mbox{ and } g_i(r,s) = \sigma_i^{-1}(r) + \frac{s-1}{c_{max}}.
\end{eqnarray*}
The interval representations of the remaining 2 interval graphs namely $I_{k+1}$ and $I_{k+2}$ 
are as follows:- 
\begin{eqnarray*}
\mbox{for every } v_{pq} \in E(G), & f_{k+1}(v_{pq}) & = [1,m]. \\
\mbox{for every } u_{pq\cdot rs} \in E(G), & f_{k+1}(u_{pq\cdot rs}) & = [j,j], \\ 
& & \mbox{where } u_{pq\cdot rs} \mbox { was obtained by } \\
& & \mbox{subdividing edge } e_j = (v_{pq},v_{rs}) \mbox{ of } H.  \\
\mbox{for every } v_{pq} \in E(G), & f_{k+2}(v_{pq}) & = [h_k(p,q),h_k(p,q)], \\ 
& & \mbox{where } h_k(p,q)= (\chi+1) - \sigma_k^{-1}(p) + \frac{q-1}{c_{max}}. \\
\mbox{for every } u_{pq\cdot rs} \in E(G), & f_{k+2}(u_{pq\cdot rs}) & = [h_k(p,q),h_k(r,s)] \mbox{, if } h_k(p,q) < h_k(r,s). \\
\mbox{for every } u_{pq\cdot rs} \in E(G), & f_{k+2}(u_{pq\cdot rs}) & = [h_k(r,s),h_k(p,q)] \mbox{, if } h_k(r,s) < h_k(p,q), \\ 
& & \mbox{where }h_k(p,q) = (\chi+1) - \sigma_k^{-1}(p) + \frac{q-1}{c_{max}} \\
& & \mbox{ and } h_k(r,s) = (\chi+1) - \sigma_k^{-1}(r) + \frac{s-1}{c_{max}}.
\end{eqnarray*} 
Observe that every edge in $G$ is of the form $(u_{pq \cdot rs}, v_{pq})$ or 
$(u_{pq \cdot rs}, v_{rs})$ where $v_{pq}$ and $v_{rs}$ are vertices of $H$ 
and $u_{pq\cdot rs}$ is the vertex introduced while subdividing edge $(v_{pq},v_{rs})$. 
Any interval graph $I_i$, where $1\leq i \leq k$, is clearly a supergraph of $G$ because in 
$f_i$ the interval corresponding to $u_{pq\cdot rs}$ has its endpoints 
on the point intervals assigned to $v_{pq}$ and $v_{rs}$. The same is true 
with interval graph $I_{k+2}$. In the interval representation $f_{k+1}$ of 
$I_{k+1}$ , any vertex $v_{pq}$ is assigned an interval $[1,m]$ which overlaps 
with the interval of every other vertex. Hence all interval graphs $I_1, I_2, 
\ldots , I_{k+2}$ are supergraphs of $G$. 

In order to show that for every $(x,y) \notin E(G)$ there exists 
some interval graph $I_i$ in our collection such that $(x,y) \notin E(I_i)$, we 
consider the following cases:\\
\textbf{case 1}: $x= v_{pq}$, $y=v_{rs}$, where $v_{pq} \neq v_{rs}$. \\
As $f_1(v_{pq}) \cap f_1(v_{rs}) = \emptyset$, 
$(v_{pq},v_{rs}) \notin E(I_1)$. \\
\textbf{case 2}: $x=u_{pq\cdot rs}, y=u_{wx\cdot yz}$, where $u_{pq\cdot rs} \neq u_{wx\cdot yz}$. \\
It is easy to verify that $f_{k+1}(u_{pq\cdot rs}) \cap f_{k+1}(u_{wx\cdot yz}) = \emptyset$ and 
hence $(u_{pq\cdot rs}, u_{wx\cdot yz}) \notin E(I_{k+1})$. \\
\textbf{case 3}: $x=u_{pq \cdot rs}, y=v_{ab}$ and $a \notin \{p,r\}$. \\
Note that $p,r,a \in [\chi]$ and since $S$ is a simply 3-suitable set for $[\chi]$, 
there exists a $\sigma_i \in S$ such that $a \notin \beta(p,r,\sigma_i)$ i.e., 
$\sigma_i^{-1}(a) < \min(\sigma_i^{-1}(p),\sigma_i^{-1}(r))$ or 
$\sigma_i^{-1}(a) > \max(\sigma_i^{-1}(p),\sigma_i^{-1}(r))$. 
$f_i(v_{ab})=[g_i(a,b),g_i(a,b)]$ and $f_i(u_{pq \cdot rs})=[g_i(p,q),g_i(r,s)]$. 
Recalling that, for any $x_1 \in [\chi]$ and $x_2 \in [c_i]$,  
$g_i(x_1,x_2) = \sigma_i^{-1}(x_1) + \frac{x_2-1}{c_{max}}$ 
it is easy to verify that $f_i(v_{ab}) \cap f_i(u_{pq \cdot rs})= \emptyset$. \\
\textbf{case 4}: $x=u_{pq \cdot rs}, y=v_{ab}$ and $a \in \{p,r\}$. \\
Assume $a=p$ (proof is similar when $a=r$). 
Assume $(v_{pb},u_{pq \cdot rs}) \in E(I_i), \forall i \in \{1, 2, \ldots, k+2\}$.
It means $(v_{pb},u_{pq \cdot rs}) \in E(I_k) \implies 
\sigma_k^{-1}(p) + \frac{q-1}{c_{max}} <  \sigma_k^{-1}(p) + \frac{b-1}{c_{max}} 
< \sigma_k^{-1}(r) + \frac{s-1}{c_{max}}  
\implies q<b$ (here we assume that $\sigma_i^{-1}(p)
< \sigma_i^{-1}(r)$. Proof is similar when $\sigma_i^{-1}(p)
> \sigma_i^{-1}(r)$). In $f_{k+2}$, note that $u_{pq \cdot rs}$ is assigned the interval 
$[(\chi+1) - \sigma_k^{-1}(r) + \frac{s-1}{c_{max}}, (\chi+1) - \sigma_k^{-1}(p) + \frac{q-1}{c_{max}}]$
and  $v_{ab}$ ($= v_{pb}$) is assigned the interval 
$[(\chi+1) - \sigma_k^{-1}(p) + \frac{b-1}{c_{max}}, (\chi+1) - \sigma_k^{-1}(p) + \frac{b-1}{c_{max}}]$. 
Therefore, $(v_{pb},u_{pq \cdot rs}) \in E(I_{k+2}) \implies b<q$. But this contradicts 
our earlier inference that $q<b$. Therefore, either $(v_{ab},u_{pq \cdot rs}) \notin E(I_k)$ 
or $(v_{ab},u_{pq \cdot rs}) \notin E(I_{k+2})$. 

We have thus shown that for any $(x,y) \notin E(G)$, $\exists i \in [k+2]$ 
such that $(x,y) \notin E(I_i)$. As each $I_i$ is a supergraph of $G$, we 
have $G=\bigcap_{i=1}^{k+2}I_i$. Applying Lemma \ref{Robertslemma}, we get 
$\boxi(G) \leq \lceil \log_2 \log_2 \chi \rceil + 3$. 
\qed
\end{proof}
\begin{corollary}
\label{sparsegraphcorollary}
Given a graph $H$, let $G$ be the graph obtained by fully subdividing $H$. Then, 
$\boxi(G) \leq \lceil \log_2 \log_2 (\Delta(H)) \rceil + 3 
 \leq \lceil \log_2 \log_2 (\Delta(G)) \rceil + 3$
\end{corollary}
\begin{proof}
By Brooks' theorem (see chapter 5 in \cite{Diest}),  
$\chi \leq \Delta(H)$ unless the graph $H$ is isomorphic to 
a complete graph $K_{\Delta(H) + 1}$ or to an odd cycle. 
If $H$ is isomorphic to $K_{\Delta(H) +1}$, then by Theorem \ref{completegraphtheorem},  
$\boxi(G) \leq \lceil \log_2 \log_2 (\Delta(H) + 1) \rceil + 2 \leq \lceil \log_2 \log_2 (\Delta(H))
 \rceil + 3$. If $H$ is an odd cycle, then $G$ will be a cycle and hence  
$\boxi(G) \leq 2 < \lceil \log_2 \log_2 (\Delta(H)) \rceil + 3$. Therefore applying  
Theorem \ref{sparsegraphtheorem}, we have $\boxi(G) \leq 
\lceil \log_2 \log_2 (\Delta(H)) \rceil + 3$. As $\Delta(H) \leq \Delta(G)$, the corollary 
follows. 
\qed
\end{proof}
\section{Line Graphs}
\label{linegraphsection}
For any bipartite graph $G$ with bipartition $\{A,B\}$, we use $C_A(G)$ 
to denote the graph with  $V(C_A(G)) = V(G)$ and 
$E(C_A(G))=E(G) \cup \{(x,y)~|~x,y \in A\}$. Thus $C_A(G)$ is 
the graph obtained from $G$ by making $A$ a clique. Similarly one can 
define $C_B(G)$. 
\begin{lemma}
\label{bipartitelemma}
For any bipartite graph $G$ with bipartition $\{A,B\}$, 
$\boxi(C_A(G)) \leq 2\cdot \boxi(G)$. 
\end{lemma}
\begin{proof}
Proof of this lemma is similar to the proof of Lemma 7 in \cite{SunMatRog}. In 
\cite{SunMatRog}  
it is proved that $\boxi(C_{AB}(G)) \leq 2\cdot \boxi(G)$, where 
$C_{AB}(G)$ refers to the graph obtained by making both $A$ and $B$ cliques. For 
the sake of completeness, we give a proof to our lemma below.  

Let $\boxi(G) = b$. Then by Lemma \ref{Robertslemma}, 
there exist $b$ interval graphs, say $I_1, I_2,$ $\ldots ,I_b$, such that 
$G = \bigcap_{i=1}^bI_i$. Let $f_i$ denote an interval representation of $I_i$, where 
$i \in [b]$. Let $s_i=\min_{x\in A}(l(f_i(x)))$ and $t_i = 
\max_{x\in A}(r(f_i(x)))$. From these $b$ interval graphs we construct $2b$ interval graphs namely 
$I_1', I_2', \ldots I_b',I_1'', I_2'', \ldots I_b''$ as follows. Let $f_i'$, $f_i''$ 
denote interval representations of $I_i'$ and $I_i''$ respectively, where $i\in [b]$. 
\begin{eqnarray*}
\mbox{Construction of } f_i': \\
\forall x\in A,~f_i'(x)& = & [s_i,r(f_i(x))]. \\
\forall x\in B,~f_i'(x)& = & f_i(x). \\
\mbox{Construction of } f_i'': \\
\forall x\in A,~f_i''(x)& = & [l(f_i(x)), t_i]. \\
\forall x\in B,~f_i''(x)& = & f_i(x). \\
\end{eqnarray*}
We claim that $C_A(G) = \bigcap_{i=1}^b (I_i' \cap I_i'')$. 
Consider any $(x,y)\in E(C_A(G))$. To show that 
$(x,y)\in E(I_i')$ and $(x,y)\in E(I_i'')$, $\forall i \in [b]$, we consider the 
following 2 cases. If $(x,y)\in E(G)$, clearly $(x,y)\in E(I_i)$. From the 
construction of $f_i'$ and $f_i''$, 
it is easy to see that $I_i'$ and $I_i''$ are supergraphs of $I_i$. 
Otherwise if $(x,y)\notin E(G)$, then $x,y \in A$ and therefore   
$[s_i,s_i] \subseteq f_i'(x) \cap f_i'(y)$ and 
$[t_i,t_i] \subseteq f_i''(x) \cap f_i''(y)$.  
%

Now, consider any $(x,y) \notin E(C_A(G)$. 
We know that $(x,y) \notin E(C_A(G)) 
\implies (x,y) \notin E(G) \implies 
(x,y) \notin E(I_i), \mbox{ for some } i\in [b]$. 
It is then easy to verify that,

(a) if $x \in A$, $y \in B$, then  $\left(f_i'(x) \cap f_i'(y) = \emptyset\right)$ or 
$\left(f_i''(x) \cap f_i''(y) = \emptyset \right)$.

(b) if $x,y \in B$, then $\left(f_i'(x) \cap f_i'(y) = \emptyset \right)$ and  
$\left(f_i''(x) \cap f_i''(y) = \emptyset \right)$.  \\
Thus we prove the claim that $C_A(G) = \bigcap_{i=1}^b(I_i' \cap I_i'')$.Therefore 
by Lemma \ref{Robertslemma}, $\boxi(C_A(G)) \leq 2 \cdot \boxi(G)$. 
\qed
\end{proof}

\begin{lemma}
\label{linegraphlemma}
Let $G$ be a bipartite graph with bipartition $\{X,Y\}$ having the following two properties:
(i) for any $y \in Y$, $d_G(y) \leq 2$ and 
(ii) for any $y_1,y_2 \in Y$, if $y_1 \neq y_2$ then $N_G(y_1) \neq N_G(y_2)$. 
Then, $\boxi(G) \leq  \lceil \log_2 \log_2 (\Delta(G))\rceil + 3$. 
\end{lemma}
\begin{proof}
If $\Delta(G) = 1$, then $G$ is a collection of isolated edges and therefore  
$\boxi(G) = 1 \leq \lceil \log_2 \log_2 (\Delta(G))\rceil + 3$. 
Let $\Delta(G) \geq 2$.    
From $G$, we construct a bipartite graph 
$G'$ with bipartition $\{X',Y'\}$
in the following way: To start with, let $G' = G$. 
For each vertex $u \in Y'$ with $d_{G'}(u) = 1$,
we add a new vertex $n_u$ to $X'$ 
such that $u$ is the only neighbour of $n_u$. For each $v \in 
Y'$ with $d_{G'}(v) = 0$,
delete $v$ from  $Y'$. 
So $X' = X \cup \{n_u~|~u \in Y \mbox{ and } d_{G}(u)=1\}$ 
and $Y' = Y \setminus 
\{v \in Y~|~v \mbox{ is an isolated vertex}\}$. We claim that $\boxi(G) \leq \boxi(G')$. 
This is because the graph obtained by removing isolated vertices from $G$ is an 
induced subgraph of $G'$ and therefore its boxicity is at most that of $G'$. As adding 
isolated vertices to any graph does not increase its boxicity, our claim follows.   

From the construction of $G'$ we can say that, for every $y \in Y'$, 
$d_{G'}(y) = 2$. 
Let $G''$ be the subgraph induced on vertices of $X'$ in ${G'}^2$, 
where ${G'}^2$ denotes the square of graph $G'$.  
It is easy to see that $G'$ can be obtained by fully subdividing 
$G''$ (Here note that if $G$ and thereby $G'$ had not satisfied property (ii), then 
the graph obtained by fully subdividing $G''$ would have just been a subgraph of $G'$). 
Therefore by our above claim and applying Corollary \ref{sparsegraphcorollary}, we get  
\[\boxi(G) \leq \boxi(G') \leq \lceil \log_2 \log_2 (\Delta(G'))\rceil + 3.\]
From the construction of $G'$ and recalling that $\Delta(G) \geq 2$, 
we infer that $\Delta(G') \leq \Delta(G)$. Therefore, 
\[\boxi(G) \leq \lceil \log_2 \log_2 (\Delta(G))\rceil + 3.\] 
\qed
\end{proof}

A \emph{critical clique} of a graph $G$ is a clique $K$ where the 
vertices of $K$ all have the same set of neighbours in $G\setminus K$, and 
$K$ is maximal under this property. Let $\mathcal{K}$ denote the collection of 
critical cliques in $G$.  The \emph{critical clique graph} of a 
graph $G$, denoted by $CC(G)$, has $V(CC(G)) = \mathcal{K}$ and 
$E(CC(G)) = \{(K_1, K_2)~|~K_1,K_2 \in \mathcal{K} \mbox{ and } 
V(K_1) \cup V(K_2) \mbox{ induces a clique in } G\}$. 
Notice that 
$CC(G)$ is isomorphic to some induced subgraph of $G$. For example, we 
can take a representative vertex from each critical clique and the induced 
subgraph on this set of vertices is isomorphic to $CC(G)$. The following 
lemma is due to Chandran, Francis and Mathew \cite{SunMatRog2} : 
\begin{lemma}
\label{criticalcliquelemma}
For any graph $G$, $\boxi(G)=\boxi(CC(G))$.
\end{lemma}

We now prove the main result of the paper. Recall that, given a 
multigraph $H$, we define its line graph $L(H)$ in the 
following way: $V(L(H)) := E(H)$ and $E(L(H)) := 
\{(e_1,e_2)~|~e_1,e_2 \in E(H),~e_1 \mbox{ and } e_2 \mbox{ share an 
endpoint in } H \}$. A graph $G$ is a line graph if and 
only if there exists a multigraph $H$ such that $G$ is isomorphic to $L(H)$. 
\begin{theorem}
\label{linegraphtheorem}
Given a multigraph $H$, let $G$ be a graph 
isomorphic to $L(H)$. Let $\Delta$ denote $\Delta(G)$ and $\chi$ represent $\chi(G)$. 
Then, 
$\boxi(G) \leq 2\Delta(\lceil \log_2 \log_2 \Delta \rceil + 3) + 1$.
\end{theorem}
\begin{proof}
Given a vertex colouring of $G$ using $\chi$ colours, 
let $D_1,D_2,\ldots ,D_{\chi}$ be the colour classes. 
For any $1\leq i \leq (\chi -1)$, let $G_i$, with 
$V(G_i)=V(G)$ and $E(G_i)=E(G) \cup \{(x,y)~|~x,y 
\in \overline{D_i}\}$, be the split graph where $D_i$ is an 
independent set and $\overline{D_i}$ a clique 
(here $\overline{D_i}=\{x\in V(G)~|~x \notin D_i\}$). 
Let $G_{\chi}^+$ be the graph 
having $V(G_{\chi}^+) = V(G)$ and 
$E(G_{\chi}^+) = \{(x,y)~|~x \in \overline{D_{\chi}}, y \in V(G)\}$. 
It is easy to see that 
\[G = G_1 \cap G_2 \cap \cdots \cap 
G_{(\chi -1)} \cap G_{\chi}^+.\]
Therefore by Lemma \ref{Robertslemma}, 
\begin{eqnarray*}
\boxi(G) & \leq & \Sigma_{i=1}^{(\chi -1)}\boxi(G_i) + \boxi(G_{\chi}^+).
\end{eqnarray*}
By Lemma \ref{criticalcliquelemma}, we know that $\boxi(G_i) = 
\boxi(CC(G_i))$. Also, observe that $G_{\chi}^+$ 
is an interval graph and hence its boxicity is 1. Therefore, 
\begin{eqnarray}
\label{criticalcliquebound}
\boxi(G) & \leq & \Sigma_{i=1}^{(\chi -1)}\boxi(CC(G_i)) + 1.
\end{eqnarray} 
We know that, $\forall i \in [(\chi -1)]$, $G_i$ is a split graph, where $D_i$ is an independent set 
and $\overline{D_i}$ a clique.  
As $CC(G_i)$ is isomorphic to some subgraph of 
$G_i$, it is also a split graph with $V(CC(G_i)) = X_i \uplus Y_i$, 
where $X_i \subseteq D_i$ is an 
independent set and $Y_i \subseteq \overline{D_i}$ a clique.   
Let $H_i$ be the bipartite graph 
obtained from $CC(G_i)$ by making $Y_i$ an 
independent set. By Lemma \ref{bipartitelemma}, we have 
$\boxi(CC(G_i)) \leq 2 \cdot \boxi(H_i)$. 
Applying this to inequality (\ref{criticalcliquebound}), we get 
\begin{eqnarray}
\label{boxicitysuminequality}
\boxi(G) \leq 2\Sigma_{i=1}^{(\chi -1)}\boxi(H_i) + 1. 
\end{eqnarray}
%
%
%
\begin{claim}
\label{linegraphclaim}
For any $i \in [(\chi -1)]$ and $y \in Y_i$, $d_{H_i}(y) \leq 2$.  
\end{claim}
\begin{proof}
Recall that $G=L(H)$ and therefore a proper vertex colouring of $G$ 
is equivalent to a proper edge colouring of $H$. 
Since in any edge colouring of $H$ a given edge $e$ cannot have 
more than 2 monochromatic neighbours, for any $y \in \overline{D_i}$, 
$|N_{G}(y) \cap D_i| \leq 2$.
Observe that the bipartite graph $H_i$ is a subgraph of $G$. Therefore,  
for any $y \in Y_i \subseteq \overline{D_i}$, 
we get $|N_{H_i}(y) \cap X_i| = |N_{H_i}(y)| = d_{H_i}(y) \leq 2$.  
\end{proof}

For any $i \in [(\chi -1)]$, $H_i$ is a bipartite graph with bipartition $\{X_i,Y_i\}$ satisfying the following 
two properties: \\
(i) by Claim \ref{linegraphclaim}, for any $ y \in Y_i$, $d_{H_i}(y) \leq 2$. \\
(ii) for any $y_1,y_2 \in Y_i$, if $y_1 \neq y_2$ then $N_{H_i}(y_1) \neq N_{H_i}(y_2)$. 
Assume for contradiction that there exist some $y_1,y_2 \in Y_i$ with $y_1 \neq y_2$ 
and $N_{H_i}(y_1) = N_{H_i}(y_2)$. Then we have $N_{CC(G_i)}(y_1) = N_{CC(G_i)}(y_2)$ 
which contradicts the fact that $CC(G_i)$ is the critical clique graph of $G_i$. 
%

Therefore by Lemma \ref{linegraphlemma}, we get 
$\boxi(H_i) \leq  \lceil \log_2 \log_2 (\Delta(H_i))\rceil + 3$. 
Since $H_i$ is a subgraph of $G$, $\Delta(H_i) \leq \Delta$. Hence, 
$$\boxi(H_i) \leq  \lceil \log_2 \log_2 \Delta\rceil + 3.$$
We thus rewrite inequality (\ref{boxicitysuminequality}) as, 
\begin{eqnarray*}
\boxi(G) \leq 2(\chi -1)(\lceil \log_2 \log_2 \Delta \rceil + 3) + 1 \leq 
2\Delta(\lceil \log_2 \log_2 \Delta \rceil + 3) + 1.
\end{eqnarray*}
As $G=L(H)$, $\Delta \leq 2(\Delta(H)-1) \leq 2(\chi -1)$. 
Therefore, 
\begin{eqnarray*}
\boxi(G) & \leq & 2(\chi -1)(\lceil \log_2 \log_2 (2(\chi -1)) \rceil + 3) + 1.
\end{eqnarray*}
\qed
\end{proof}
%
%
\section{Lower Bound for Boxicity of a Hypercube}
For any non-negative integer $d$, a $d$-dimensional 
hypercube $H_d$ has its vertices corresponding to 
the $2^d$ binary strings each of length $d$. 
Two vertices are adjacent if and only if their binary strings 
differ from each other in exactly one bit position. 
\begin{theorem}
$\boxi(H_d) \geq \frac{\lceil \log_2 \log_2 d \rceil + 1}{2}$
\end{theorem}
\begin{proof}
For any vertex $v \in V(H_d)$, 
let $g(v)$ denote the number of ones in the bit 
string associated with $v$. Let $X=\{v \in V(H)~|~g(v)=1 \mbox{ or }
g(v)=2\}$. Let $H'$ be the subgraph of $H$ induced on the vertex set 
$X$. 
We can see that $H'$ is a bipartite graph with bipartition 
$\{A, B\}$, where $A = \{v\in V(H')~|~g(v)=1\}$ and 
$B=\{v\in V(H')~|~g(v)=2\}$. 


It is easy to observe that $H'$ is a graph obtained by fully subdividing 
$K_{|A|}$, where $K_{|A|}$ refers to a complete graph on $|A|=d$ vertices. 
Then by Theorem \ref{completegraphtheorem}, we can say that 
\[\boxi(H') \geq \frac{\lceil \log_2 \log_2 d \rceil + 1}{2}.\] 
As $H'$ is an induced subgraph of $H$, 
\[\boxi(H) \geq \boxi(H') \geq \frac{\lceil \log_2 \log_2 d \rceil + 1}{2}.\] 
\qed
\end{proof}
\bibliography{mathewref}
\end{document}